\documentclass{amsproc}
\usepackage{amsmath}
\usepackage{amsfonts}
\usepackage{latexsym}
\usepackage{amssymb}
\usepackage{diagrams}
\usepackage[dvips]{graphicx}

\newcommand{\bdm}{\begin{displaymath}}
\newcommand{\edm}{\end{displaymath}}
\newcommand{\R}{\mathbb{R}}
\newcommand{\Z}{\mathbb{Z}}
\newcommand{\co}{\colon\thinspace}

\newcommand{\TC}{\mathbf{TC}}
\newcommand{\TCS}{\mathbf{TC}^S}
\newcommand{\EG}{EG\times_G}
\newcommand{\wgt}{\mathrm{wgt}}
\newcommand{\supp}{\rm {supp}}
\newcommand{\cl}{\rm {cl}}

\newtheorem{theorem}{Theorem}
\newtheorem{lemma}[theorem]{Lemma}

\theoremstyle{definition}
\newtheorem{definition}[theorem]{Definition}
\newtheorem{example}[theorem]{Example}

\newtheorem{proposition}[theorem]{Proposition}
\newtheorem{corollary}[theorem]{Corollary}

\theoremstyle{remark}
\newtheorem{remark}[theorem]{Remark}

\begin{document}

\title{Symmetric Motion Planning}

\author{Michael Farber}
\address{Department of Mathematical Sciences, University of Durham, Durham DH1 3LE}

\email{Michael.Farber@durham.ac.uk}
\thanks{The first author was supported by a grant from the Royal Society.
The second author was supported by a grant from the UK
  Engineering and Physical Sciences Research Council.}

\author{Mark Grant}
\address{Department of Mathematical Sciences, University of Durham, Durham DH1 3LE}
\email{Mark.Grant@durham.ac.uk}

\date{January 1, 1994 and, in revised form, June 22, 1994.}



\begin{abstract}
In this paper we study symmetric motion planning algorithms, i.e.
such that the motion from one state $A$ to another $B$, prescribed
by the algorithm, is the time reverse of the motion from $B$ to
$A$. We experiment with several different notions of topological
complexity of such algorithms and compare them with each other and
with the usual (non-symmetric) concept of topological complexity.
Using equivariant cohomology and the theory of Schwarz genus we
obtain cohomological lower bounds for symmetric topological
complexity. One of our main results states that in the case of
aspherical manifolds the complexity of symmetric motion planning
algorithms with fixed midpoint map exceeds twice the cup-length.

We introduce a new concept, the sectional category weight of a
cohomology class, which generalises the notion of category weight
developed earlier by E. Fadell and S. Husseini. We apply this
notion to study the symmetric topological complexity of aspherical
manifolds.
\end{abstract}

\maketitle

\section{Introduction}

The Motion Planning Problem is a central theme in Robotics, which
invites applications of tools of algebraic topology. Any
autonomous mechanical system or robot capable of performing tasks
must first be told how to move between different states of the
system. A {\em Motion Planner} in a given mechanical system is a
rule which assigns to any pair of states $A$ and $B$ of the system
a continuous motion from $A$ to $B$. Topology enters the picture
by regarding the set of admissible states of a mechanical system
as the points of a topological space, $X$, called the {\em
  configuration space} of the system. In practice this space is
usually determined by several real parameters, and so has the
natural topology induced by some Euclidean metric. The system
becomes synonymous with the space $X$, and continuous motions of
the system are represented by continuous paths in $X$ (continuous
maps of the unit interval $I=[0,1]$ to $X$).

Let $PX$ denote the space of all continuous paths in $X$, endowed
with a suitable topology (such as the compact-open topology). The
map $\pi\co PX\to X\times X$ which takes a path in $X$ to its
end-points, given by $\pi(\gamma)=(\gamma(0),\gamma(1))$, is a
fibration in the sense of Serre. A Motion Planner in $X$ is then a
function $s$ from the Cartesian product $X\times X$ to the path
space $PX$, such that the composition $\pi\circ s$ is the identity
function on $X\times X$. Whilst such discontinuous functions
always exist when $X$ is path-connected, one may show that a {\em
continuous} Motion Planner in $X$ (which is therefore a continuous
section of the fibration $\pi$) exists if and only if $X$ is
contractible \cite{Far03}. Thus Motion Planners in $X$ may have
essential discontinuities, which reflect the homotopy properties
of $X$ and provide a measure of the complexity of the task of
navigation in $X$.

 The problem has been studied extensively from this viewpoint
 by the first author in \cite{Far03}, \cite{Far04}, and \cite{Far06},
 where a new homotopy invariant was explored.
 For any space $X$, the {\em Topological Complexity} of $X$ is a positive natural
 number $\TC(X)$, which may be defined in a number of equivalent ways.
 The definition we will give is based on the following concept due to
 A.\ S.\ Schwarz \cite{Sch66}.
 \begin{definition}
 Let $p\co E\to B$ be a continuous map and $U\subseteq B$. Recall that a map $s\co
U \to E$ is a section of $p$ over $U$ if $p\circ s=\mathbf{1}_U$,
the identity on $U$. The {\em Schwarz genus} of $p$, denoted
$g(p)$, is the minimum cardinality
 among coverings of $B$ by open sets, over each of which $p$ has a continuous
 section (in \cite{CLOT} the term {\it sectional category} was used).
 \end{definition}

\begin{definition}
The {\em Topological Complexity} of a path-connected topological
space $X$, denoted $\TC(X)$, is the Schwarz genus of the path
fibration $\pi\co PX\to X\times X$.
\end{definition}
Thus $\TC(X)=g(\pi)\leq k$ iff there is an open cover $\{
U_1,\ldots, U_k\}$ of $X\times X$ and continuous maps $s_i\co
U_i\to PX$ such that $\pi\circ s_i=\mathbf{1}_{U_i}$, for
$i=1,\ldots ,k$. The number $\TC(X)$ depends only on the homotopy
type of $X$ (see \cite{Far03}) and provides a measure of the
intrinsic complexity of the Motion Planning problem in $X$. This
invariant is similar in spirit to the Lusternik-Schnirelmann
category $\mathrm{cat}(X)$, and, whilst the two are independent,
they satisfy the inequalities
 \bdm
 \mathrm{cat}(X)\leq\TC(X)\leq\mathrm{cat}(X\times X).
 \edm

 In addition, there are cohomological lower bounds for $\TC(X)$ based
 on `zero-divisors cup-length' (see \cite{Far03} for details) and one
 may obtain upper bounds, for instance by constructing explicit Motion
 Planners in $X$. With these considerations one may calculate the
 Topological Complexity of a large number of spaces, including
 spheres, simply connected symplectic manifolds (such as
 $\mathbb{C}P^n$, $n\geq 1$) and tori. The paper \cite{Far06} contains
 a recent survey of such results, and along with \cite{Sch66} is our
 basic reference.

In this paper we study several variations on the Motion Planning
problem. We impose additional, quite natural, symmetry constraints
on our Motion Planners, namely that the motion from $A$ to $A$
should be constant at $A$, while the motion from $B$ to $A$ should
be the motion from $A$ to $B$ traversed in the opposite direction.
This is formalised as follows.
\begin{definition}\label{SMP}
A {\em Symmetric Motion Planner} in $X$ is a (possibly
discontinuous) function $s\co X\times X\to PX$ such that $\pi\circ
s=\mathbf{1}_{X\times X}$ and the following conditions are
satisfied for all $t\in I$: \begin{eqnarray}\label{symmetric}
s(A,A)(t)=A,\quad s(B,A)(t)=s(A,B)(1-t).
\end{eqnarray}
\end{definition}
Armed with this Definition, we will define in Section 2 a new
invariant which measures the intrinsic complexity of Symmetric
Motion Planning in a topological space. The {\em Symmetric
Topological Complexity} of $X$, denoted $\TCS(X)$, is compared
with the usual $\TC(X)$ and some examples are explored. Note that
$\TCS(X)$, unlike $\TC(X)$, is not homotopy invariant, and harder
to deal with. The computation of this new invariant requires
usable cohomological lower bounds like the one mentioned above for
$\TC(X)$; this is done in Section \ref{cohbx} of the present
paper, using $\mathbb{Z}_2$-equivariant cohomology and results of N.\ E.\ Steenrod and A.\ Haefliger. We find that for a closed smooth manifold $X$,
$$\TCS(X)\geq \textrm{\rm {cup-length}}(N)+2,
$$ where $N$ denotes the sub-ring of $H^*(X)\otimes H^*(X)$ spanned by
the norm elements (elements $x\otimes y+y\otimes x$ with $x\neq y$).

In Section \ref{midpointsec} we impose further constraints on our
Symmetric Motion Planners, namely that the mid-point of a motion
between $A$ and $B$ should depend continuously on $A$ and $B$.
This has the effect of greatly increasing the complexity of
navigation in certain cases. In particular we consider Symmetric
Motion Planners such that the mid-point of each motion is a fixed
base state $A_0\in X$. One of our main results, Theorem
\ref{cuplength}, asserts that the Symmetric Topological Complexity
with constant mid-point of a closed aspherical manifold $X$
exceeds twice the cup-length of $X$. We find that in the case of
planar robot arm with $n$ revolving joints (having the torus $T^n$
as the configuration space) the symmetric topological complexity
is $2n+1$ while $\TC(X)=n+1$; thus in this instance the
requirement that motions are symmetric and return to a fixed base
state increases the complexity of navigation by $n$.


In Section  \ref{secsec} (which can be read independently of the
rest of the paper) we employ a new tool for estimating the Schwarz
genus of a fibration. We introduce the notion of sectional
category weight of a cohomology class. It is a natural
generalization of the usual category weight defined by Fadell and
Husseini \cite{FH92} and developed by Rudyak \cite{Rud99}. One may
improve on the classical cohomological lower bound for the genus,
given in \cite{Sch66}, by finding cohomology classes of weight at
least two.

In Section \ref{secaspherical} we show that some cohomology
classes which arise in the $\Z_2$-equivariant cohomology have
sectional category weight $\geq 2$. This is the main ingredient of
the proof of Theorem \ref{cuplength}.

{\em Everywhere in this paper we denote by $G$ the cyclic group
$\mathbb{Z}_2$ of order two. Cohomology is taken with coefficients in
$G$ unless otherwise stated.}

\section{Symmetric Topological Complexity}
Let $X$ be a path-connected polyhedron. The path fibration
\begin{eqnarray}\label{pi} \pi\co PX\to X\times X\end{eqnarray}
restricts to a fibration
\begin{eqnarray}\label{piprime}
\pi'\co
P'X\to F(X;2),\end{eqnarray}
where $F(X;2)=\{(x,y)\in X\times
X\mid x\neq y\}$ is the space of ordered pairs of distinct points
in $X$, and $P'X$ is the subspace $\{\gamma\co I\to
X\mid\gamma(0)\neq\gamma(1)\}\subseteq PX$ consisting of paths
with distinct endpoints.

The spaces $P'X$ and $F(X;2)$ carry free $G$-actions, defined in
the latter case by permutation of factors and in the former by
sending a path $\gamma$ to its inverse $\overline{\gamma}$, given
by $\overline{\gamma}(t)=\gamma(1-t)$. Note that $\pi'\co P'X\to
F(X;2)$ is an equivariant map of free $G$-spaces. Hence the
quotient map \begin{eqnarray}\label{map} \pi_G:=\pi'/G\co P'X/G\to
B(X;2)
\end{eqnarray} is also a fibration, where $B(X;2)$ denotes the orbit space
$F(X;2)/G$ of unordered pairs of distinct points in $X$.

Comparing with Definition \ref{SMP}, we see that a Symmetric
Motion Planner in $X$ describes a function $s\co B(X;2)\to P'X/G$
such that $\pi_G \circ s$ is the identity map on $B(X;2)$.
Conversely, such a function $s$ completely describes a Symmetric
Motion Planner, since the latter must map a point $(A,A)$ on the
diagonal of $X\times X$ to the constant path at $A$. Of course a
{\em continuous} map $s$ of this kind may exist only in very few
cases. We wish to measure the essential discontinuities of
Symmetric Motion Planning in $X$. Thus we are led to the following
definition.
\begin{definition}\label{def1}
The {\em Symmetric Topological Complexity} of $X$, denoted
$\TCS(X)$, is defined to be one plus the Schwarz genus of the
fibration (\ref{map}). In other words,
\begin{eqnarray}\TCS(X) = 1+g(\pi_G).\end{eqnarray}
\end{definition}

We adopt the convention that the Schwarz genus of $p\co E\to B$
vanishes iff $E=\emptyset=B$. Note that $B(X,2)$ is empty if and
only if $X$ is a single point; in this case $\TCS(X)=1$. If $X$ is
not a single point then $g(\pi_G)\geq 1$ and therefore
\begin{eqnarray}\TCS(X)\geq 2.
\end{eqnarray}

\begin{example} Let $X\subseteq \R^n$ be a convex subset, not a single point. For $A,
B\in X$ define $s(A,B)(t) = (1-t)A+tB$ where $t\in I=[0,1]$. This
defines a continuous equivariant section of (\ref{map}) and hence
$\TCS(X)=2$.
\end{example}

\begin{example}
Let $X$ be a finite tree. We may view $X$ as a metric space by
specifying length of every edge. Then for any pair of points $A,
B$ there is a unique constant speed curve of minimal length $s(A,
B)\co [0,1]\to X$ starting at $A$ and ending at $B$. We obtain a
continuous equivariant section of (\ref{map}) which implies
$\TCS(X)=2$.
\end{example}

\begin{example} More generally, it is easy to see that
for any contractible $X$ which is not a single point one has
$\TCS(X)=2$. Indeed, since $X$ is contractible there exists a
continuous map $x\mapsto \gamma_x\in PX$ such that $\gamma_x(0)=x$
and $\gamma_x(1)=x_0$. Then setting $s(A,B)$ to be equal the
concatenation of $\gamma_A$ and the inverse path to $\gamma_B$
gives a symmetric equivariant section of (\ref{piprime}).

One expects that the converse statement is true as well: any
path-connected polyhedron $X$ with $\TCS(X)=2$ is contractible.
\end{example}

The following lemma will be used to justify Definition \ref{def1}.

\begin{lemma}\label{lm1} Let $U$ be an open subset $U\subseteq B(X,2)$. Denote
by $\tilde U$ the preimage of $U$ under the projection $q\co
F(X,2)\to B(X,2)$. Any continuous section $s\co U\to P'X/G$ of
fibration (\ref{map}) over $U$ determines a continuous equivariant
section $\tilde s\co \tilde U\to P'X$ of fibration (\ref{piprime})
over $\tilde U$.
\end{lemma}
\begin{proof}
Let $(A,B)\in \tilde U$. Then $s(q(A,B))\in P'X/G$ is an equivalence
class containing two paths $\gamma$ and $\overline\gamma$, one of
which goes from $A$ to $B$ and the other of which goes from $B$ to
$A$. Without loss of generality we may assume that $\gamma(0)=A$ and
$\gamma(1)=B$, and set $\tilde{s}(A,B)=\gamma$. In this way one
describes a continuous equivariant section $\tilde s\co \tilde U\to
P'X$.
\end{proof}

We now show that if $k= \TCS(X)-1$ then we may partition $X\times X$ into $k$
{\em disjoint} subsets with a continuous
symmetric motion planner on each. Given an open cover $\{U_1, \dots, U_k\}$
of $B(X,2)$ and a sequence of continuous section $s_i\co U_i\to
P'X/G$, one uses Lemma \ref{lm1} to obtain an open cover $\tilde
U_1, \dots, \tilde U_k$ of $F(X,2)$ with equivariant continuous
sections $\tilde s_i\co \tilde U_i\to P'X$, where $i=1, \dots, k$.
Consider a partition of unity $\{f_1, \dots, f_k\}$ subordinate to
the cover $\{U_1, \dots, U_k\}$. Here $f_i\co B(X,2)\to \R_+$ is a
continuous function with $\supp(f_i)\subseteq U_i$ and
$f_1+\dots+f_k=1$. Composing with the projection $q\co F(X, 2)\to
B(X,2)$ gives an equivariant partition $\tilde f_i = f_i\circ q$
of unity of $F(X,2)$, where $i=1, \dots, k$, subordinate to
$\{\tilde U_1, \dots, \tilde U_k\}$. Define sets $V_i\subseteq
F(X,2)$ by
\begin{eqnarray*}
(A,B)\in V_i & \iff &\left\{
\begin{array}{lll}
\tilde f_i(A,B) \geq 1/k,\\
\\
\tilde f_j(A,B) < 1/k &\mbox{for all}& j<i .
\end{array}
\right.
\end{eqnarray*}
Then

(a) each $V_i$ is involution invariant and is contained in $\tilde
U_i$;

(b) the sets $V_i$ are pairwise disjoint;

(c) $V_1\cup \dots \cup V_k=F(X,2)$.

As a result one obtains a symmetric motion planning algorithm
$$s\co X\times X\to PX$$ by setting $s(A,B) = \tilde s_i(A,B)$ if
$A\not= B$ and $(A, B)\in V_i$; we also set $s(A,A)(t)=A$.

\begin{corollary}\label{cor1} One has
$
 \TC(X)\leq \TCS(X).
$
\end{corollary}
\begin{proof} In notations introduced in the paragraph preceding
Corollary \ref{cor1}, the sets $\tilde U_1, \dots, \tilde U_k$
constitute an open cover of $F(X,2)=X\times X - \Delta$ where
$\Delta$ is the diagonal, with continuous sections $\tilde s_i$
over each $\tilde U_i$. Hence it is enough to show that there is
an open neighbourhood $\tilde U_0\subseteq X\times X$ of $\Delta$
which supports a continuous equivariant section $\tilde s_0\co
\tilde U_0\to PX$ which takes each pair $(A,A)$ into the constant
path at $A$, for $A\in X$.

Since $X$ is an ENR, we may find an embedding $e\co X\to\R^n$ and an open
neighbourhood $e(X)\subseteq N\subseteq \R^n$ which admits a retraction
$r\co N\to X$ onto $X$. Let $\tilde U_0$ be a neighbourhood of
$\Delta$ in $X\times X$ which is involution invariant and such
that for all $(A,B)\in \tilde U_0$ the straight line segment
connecting $e(A)$ and $e(B)$ is contained in $N$. Now we define the
desired section $\tilde s_0\co \tilde U_0\to PX$ by $ \tilde s_0(A,
B)(t) = r((1-t)e(A)+te(B)). $ This completes the proof.
\end{proof}

Corollary \ref{cor1} implies that all cohomological lower bounds
for $\TC(X)$ (see \cite{Far03}) are valid also for $\TCS(X)$.

\begin{proposition}\label{dim}
For any $X$ we have $\TCS(X)\leq 2\mathrm{dim}(X)+2$ where
$\mathrm{dim}(X)$ denotes covering dimension of $X$. If $X$ is a
closed smooth manifold then $\TCS(X)\leq 2\mathrm{dim}(X)+1$.
\end{proposition}
\begin{proof}
Since the genus of a fibration may not exceed the category of the base
(see Theorem 5 of \cite{Sch66}), which in turn may not exceed the
dimension of the base plus one, we obtain
$
\TCS(X)\leq\mathrm{cat}(B(X;2))+1\leq\mathrm{dim}(B(X;2))+2\leq
2\mathrm{dim}(X)+2.
$

Let $X$ be a closed manifold with $\mathrm{dim}(X)=n$. Let
$Y\subseteq X\times X$ denote the space obtained by removing an
open $G$-invariant tubular neighbourhood of the diagonal. Then $Y$
inherits a free $G$-action, and is $G$-equivariantly homotopy
equivalent to $F(X;2)$. Its quotient $Y/G=Y'$ is therefore
homotopy equivalent to $B(X;2)$, and is a compact $2n$-manifold
with boundary. A standard Morse theoretical argument gives that
$Y'$ is homotopy equivalent to a complex $Y''$ of dimension
$2n-1$. Thus $
\TCS(X)\leq\mathrm{cat}(B(X;2))+1=\mathrm{cat}(Y'')+1\leq\mathrm{dim}(Y'')+2=2n+1.
$
\end{proof}
\begin{proposition} \label{sphere}
$\TCS(S^n)\leq 3$ for any $n\geq 1$.
\end{proposition}
\begin{proof}
The inclusion $S^n\hookrightarrow F(S^n;2)$ given by $A\mapsto (A,-A)$
is a $G$-equivariant homotopy equivalence, where $G$ acts on the
sphere antipodally. We may describe the homotopy inverse $\phi\co F(S^n;2)\to S^n$ of this map as follows. Let $d$ denote the standard metric on $S^n$.
Given any point $(A,B)\in F(S^n;2)$,
choose a geodesic great circle $\Gamma$ containing both $A$ and $B$
(note that if $A$ and $B$ are non-antipodal, so $A\neq -B$, then there
is no choice to make). Then there is a unique pair of points
$(A',B')$, both on $\Gamma$, which are antipodal ($B'=-A'$) and such
that $d(A,A')=d(B,B')<d(A,B')$. We set $\phi(A,B)=A'$, and
note that that this is independent of the choice of $\Gamma$ since if
$A=-B$ then $A=A'$. Also note that $\phi$ is equivariant, since $\phi(B,A)=B'=-A'$.
\begin{center}
\resizebox{4cm}{3.5cm}{\includegraphics[130,365][434,633]{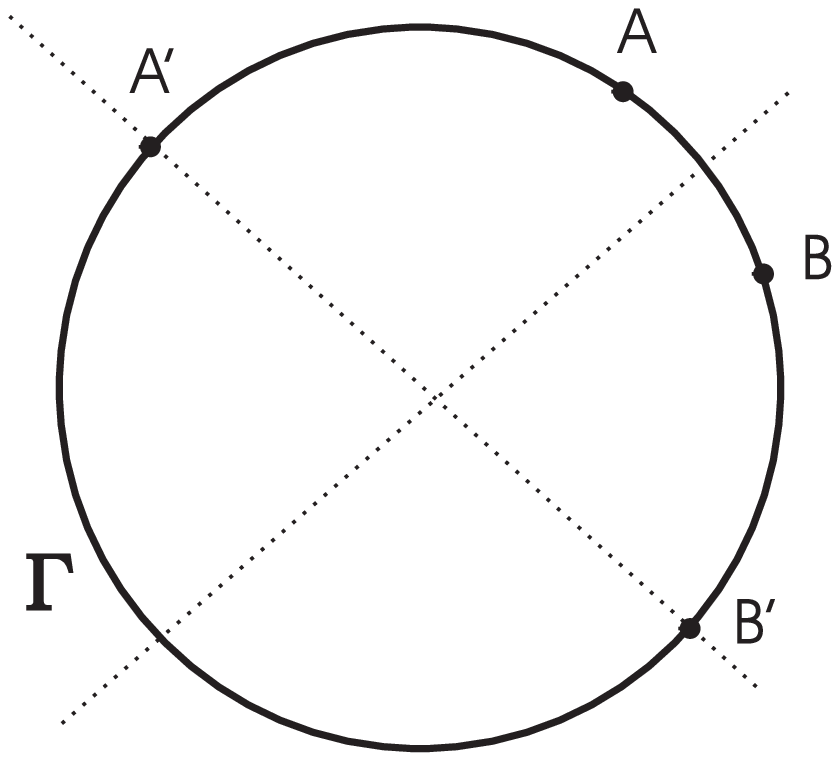}}
\end{center}

We must show that $g(\pi_G)\leq 2$. Clearly, it suffices to cover
$F(S^n;2)$ by two $G$-invariant open sets, over each of which the restricted
end-point map $\pi'\co P'S^n\to F(S^n;2)$ has an {\em equivariant}
section. Fix a point $A_0\in S^n$, and let $U\subseteq S^n$ be a small
open disk centred at $A_0$. Then $-U:=\{-C\mid C\in U\}$ is a small
open disk centred at $-A_0$. The sets
\bdm
W_1=\phi^{-1}(S^n-\{A_0,-A_0\}),\quad W_2=\phi^{-1}U\cup\phi^{-1}(-U)
\edm
are clearly open and $G$-invariant, and cover $F(S^n;2)$. We will
describe equivariant sections $s_1$, $s_2$ of $\pi'$ over these sets.

Given $(A,B)\in W_1$, choose a geodesic great circle $\Gamma$ and pair of antipodal points $(A',B')$ on $\Gamma$, as in the construction of $\phi$. Let $\alpha$ (resp.\ $\beta$) denote the path along $\Gamma$ from $A$ to $A'$ (resp.\ from $B$ to $B'$). Since $\{A',B'\}\neq\{A_0,-A_0\}$ there is a unique path $\xi$ from $A'$ to $B'$ which travels along the geodesic great circle through $A_0$. We may then set
$
s_1(A,B)=\alpha\circ\xi\circ\overline\beta.
$
This describes a continuous equivariant section $s_1\co W_1\to P'X$.

Given $(A,B)\in W_2$, the pair $(A',B')$ lies in either $U\times -U$ or
$-U\times U$. In either case, there is a unique path $\eta$ from $A'$ to $B'$
of constant speed, which travels first along the geodesic arc from $A'$
to the centre of the disk containing $A'$, then along some fixed geodesic
arc between $A_0$ and $-A_0$, then from the centre of the disk containing
$B'$ to $B'$. It is not hard to see that setting
$
s_2(A,B)=\alpha\circ\eta\circ\overline{\beta}
$
describes a continuous equivariant section $s_2\co W_2\to P'X$ of $\pi'$ over $W_2$.
\end{proof}

In the next section we shall see that $\TCS(S^n)=3$ for all $n\geq
1$. This should be compared with the corresponding result for
$\TC(S^n)$ which equals $2$ for $n$ odd and $3$ for $n$ even.

\section{Cohomology of $B(X;2)$}\label{cohbx}
The main result of this Section, Theorem \ref{lower}, gives a lower
bound for $\TCS(X)$ in terms of the structure of the cohomology
algebra of $X$ with $\Z_2$ coefficients, assuming that $X$ is a
closed smooth manifold.

We use the following result of A. Schwarz (\cite{Sch66}, Theorem
4).
\begin{theorem}\label{cohom}
Let $p\co E\to B$ be a continuous map. Suppose there exist
cohomology classes $\xi_1,\ldots ,\xi_k\in H^*(B)$ with any
coefficients, such that $p^*(\xi_1)=\ldots =p^*(\xi_k)=0$ and the
product $\xi_1\cup\cdots\cup\xi_k$ is different from zero. Then
$g(p)> k$.
\end{theorem}
We plan to apply this theorem to fibration (\ref{map}). To make
this work we need a good understanding of the cohomology algebra
of $B(X;2)$. A. Haefliger \cite{Hae61}, who studied embeddings and
immersions of manifolds, gave a very explicit description of
$H^\ast(B(X,2))$. In this Section we recall his results and apply
them to our problem of estimating the Schwarz genus of
(\ref{map}).

For the remainder of this Section $X$ will denote a closed
manifold and all cohomology will be taken with coefficients in the
group $G=\mathbb{Z}_2$ which will be skipped from the notation.
Let $\rho\co EG\to BG$ be a universal principal $G$-bundle. Hence
$EG$ is a contractible space with a free right $G$-action, and
$\rho$ is the orbit map to $BG=EG/G$. Recall that the cohomology
algebra of the group $G$ is
$H^*(G)=H^*(BG)\cong\mathbb{Z}_2[\mu]$, the polynomial algebra
over $G$ on one generator $\mu$ of degree 1.

The inclusion $F(X;2)\hookrightarrow X\times X=X^2$ is
$G$-equivariant, as are the fibrations $\pi$ and $\pi'$ (see
(\ref{pi}) and (\ref{piprime})), allowing us to form the following
diagram of homotopy orbit spaces.
\begin{eqnarray}\label{diag1}
\begin{array}{ccccc}
P'X/G & \simeq & EG\times_G P'X &\rTo & \EG PX \\ \\
\pi_G \downarrow & &\downarrow 1\times_G \pi' & &\downarrow
1\times_G\pi \\ \\ B(X;2) &\simeq & \EG F(X;2) &\rTo & \EG X^2
\end{array}
\end{eqnarray}
If $X$ and $Y$ are right and left $G$-spaces respectively, then
$X\times_G Y$ denotes the orbit space of $X\times Y$ under the
diagonal $G$-action. The homotopy equivalences on the left of the
diagram result from the fact that $F(X;2)$ and $P'X$ are free
$G$-spaces. The next Theorem describes the equivariant (Borel)
cohomology ring $H^*(\EG X^2)$ of the $G$-space $X^2$, where the
$G$-action on $X^2$ permutes the factors.

Recall that all cohomology groups are with coefficients in
$G=\Z_2$.

\begin{theorem}\label{Steen} Let $X$ be a finite
polyhedron. Then there is an isomorphism of $H^*(G)$-algebras
\begin{eqnarray}
H^*(\EG X^2)\cong H^*(G, H^*(X)\otimes H^*(X)).
\end{eqnarray}
Here we view $H^*(X)\otimes H^*(X)$ as a $G$-module, where the
non-trivial element of $G$ acts by $x\otimes y\mapsto y\otimes x$ (no
signs are introduced here since we are working mod 2). Then the right
hand side is cohomology of the group $G$ with coefficients in this
$G$-module.
\end{theorem}
 A. Haefliger credits this theorem to an unpublished work
of N.\ E.\ Steenrod. One may interpret its statement as saying
that the $G$-cohomology Cartan-Leray spectral sequence of the
regular covering $EG\times X^2\to\EG X^2$ collapses at the
$E_2$-term (see \cite{Bro82} or \cite{McC01}). Below we describe
the structure of this $E_2$-term $H^*(G, H^*(X)\otimes H^*(X))$
explicitly in terms of the cohomology ring $H^*(X)$.

Consider the ring $\left(H^*(X)\otimes H^*(X)\right)^G$ of
invariants of $H^*(X)\otimes H^*(X)$. It contains a subring
consisting of {\em diagonal elements}, which are linear
combinations of elements of the form $x\otimes x$ where $x\in
H^*(X)$; denote this subring by $D$. It also contains a subring
$N$ consisting of linear combinations of elements of the form
$x\otimes y +y\otimes x$ where $x\neq y\in H^*(X)$, the so called
{\em norm elements}. Then as $H^*(G)$-algebras, we have
\begin{eqnarray}
H^*(G,H^*(X)\otimes H^*(X))\cong H^*(G)\otimes D +N.
\end{eqnarray} This algebra is generated additively by elements of
the form $\alpha=\mu^i\otimes x\otimes x$ where $x\in H^*(X)$ (the
exponent $i\geq 0$ of $\mu$ is called the {\em $G$-degree} of
$\alpha$), and the norm elements $x\otimes y+y\otimes x\in N$
(which have $G$-degree zero). It is generated as a
$H^*(G)$-algebra by the sub-ring consisting of elements of
$G$-degree 0, which is precisely the ring of invariants
$\left(H^*(X)\otimes H^*(X)\right)^G$. The action of $H^*(G)$ is
free on $D$ and trivial on $N$.

One may now investigate the cohomology algebra $H^*(B(X;2))$ with
the aid of the exact sequence \bdm \cdots\to H^*(\EG X^2,\EG
F(X;2))\to H^*(\EG X^2)\stackrel{j^*}{\to} H^*(B(X;2))\to\cdots
\edm (Here $j\co B(X;2)\hookrightarrow \EG X^2$ is the composition
along the bottom of diagram (\ref{diag1})). Such an investigation
leads to the main theorem of \cite{Hae61}.

To describe the statement of the theorem, we must introduce
several maps. We assume below that $X$ is a closed smooth
$n$-dimensional manifold. Fixing a point $e\in EG$ we obtain
inclusions \bdm r\co X^2\hookrightarrow\EG X^2,\quad r_0\co
X\hookrightarrow BG\times X=\EG X, \edm given by $r[A,B]=[e,A,B]$
and $r_0[A]=([e],A)$, which up to homotopy are independent of the
choice of $e$. The induced maps $r^*\co H^*(\EG X^2)\to H^*(X^2)$
and $r_0^*\co H^*(BG)\otimes H^*(X)\to H^*(X)$ in cohomology are
given by \bdm r^*(x\otimes y+y\otimes x)=x\otimes y+y\otimes
x,\quad r^*(\mu^i\otimes x\otimes x)=\left\{\begin{array}{ll}
        x\otimes x &\textrm{if $i=0$,}\\
        0   &\textrm{if $i>0$,}
\end{array}\right.
\edm
\bdm
r_0^*(\mu^i\otimes x)=\left\{ \begin{array}{ll}
       x &\textrm{if $i=0$,}\\
        0   &\textrm{if $i>0$,}
\end{array}\right.
\edm and can be described as `taking the $G$-degree 0 part'. We
will also consider the {\em generalised diagonal map} \bdm
\Delta_G\co BG\times X\to \EG X^2,\quad \Delta_G([e],A)=[e,A,A].
\edm The induced map
$$\Delta_G^*\co H^*(\EG X^2)\to H^*(BG\times X)$$ plays an
important r\^{o}le in the construction of the Steenrod square
cohomology operations
$$\mathrm{Sq}^i\co H^*(X)\to H^{*+i}(X).$$ It may be completely described
by noting that it is an $H^*(G)$-module homomorphism which
vanishes on $N$, and for $x\in H^k(X)$ we have \bdm
\Delta_G^*(1\otimes x\otimes x)=\sum_{i=0}^k
\mu^{k-i}\otimes\mathrm{Sq}^i(x)\in H^{2k}(BG\times X). \edm The
diagonal embedding $\Delta\co X\hookrightarrow X^2 $ induces
a {\em
  Gysin} or {\em pushforward} map
\bdm \Delta_!\co H^{*-n}(X)\to H^*(X^2) \edm which is given by
$\Delta_!(x)=(1\times x)\cup\delta$, where $\delta\in H^n(X^2)$ is
the diagonal class (for definitions and properties of this class
see \cite{MS74}, Chapter 11).

Finally, there is the map
$$\varphi\co H^{*-n}(BG\times X)\to H^*(BG\times
X)$$ given by cup product with the element \bdm \sum_{l=0}^n
\mu^{n-l}\otimes w_l\in H^n(BG\times X), \edm where $w_l$ for
$l=0,\ldots ,n$ is the $l$-th tangential Stiefel-Whitney class of
$X$ \cite{MS74}. Note that $\phi$ is injective since $w_0=1$.

\begin{theorem}[A. Haefliger \cite{Hae61}] The following
diagram commutes and has exact rows.
\begin{eqnarray*}
\begin{array}{ccccccccc}
0&\to & H^{*-n}(X) & \stackrel {\Delta_!}\longrightarrow &
H^*(X^2) & &  & & \\ \\
   & & \uparrow {r_0^*}&  &       \uparrow{r^*}   & & & & \\ \\
0&\to & H^{*-n}(BG\times X)& \rTo & H^*(\EG X^2) &\rTo^{j^*}
   H^*(B(X;2))&\to 0 \\ \\
    & & &\searrow{\varphi} & \downarrow{\Delta_G^*} & & &  \\ \\
  & & & &  H^*(BG\times X) & & &
\end{array}
\end{eqnarray*}
An element $\alpha\in H^*(\EG X^2)$ is in the kernel of $j^*$ if
and only if: (a) There exists an element $\beta\in
H^{*-n}(BG\times X)$ such that $\varphi(\beta)=\Delta_G^*(\alpha)$
(such an element $\beta$ is unique), and (b)
$\Delta_!r_0^*(\beta)=r^*(\alpha)$.
\end{theorem}
This theorem facilitates explicit calculations of the algebra
$H^*(B(X;2))$ for certain manifolds $X$ for which the
Stiefel-Whitney classes of $X$ and the action of the Steenrod
algebra on $H^*(X)$ are known (the reader may enjoy checking for
example that $H^*(B(S^n;2))\cong\Z_2[\mu]/(\mu^{n+1})$, which
follows from the $G$-equivariant homotopy equivalence $S^n\to
F(S^n;2)$). A general description of $H^*(B(X;2))$ in terms of
these data, valid for all $X$, was obtained by Yo Ging-tzung
\cite{Ging63}. Unfortunately the details are messy and not so
instructive. Here we content ourselves with the following.
\begin{corollary}\label{Nsurvives} $H^*(B(X;2))$ contains a sub-ring $\tilde{N}$
  isomorphic to $N$.
\end{corollary}
\begin{proof}
Consider a non-zero element $\alpha=x\otimes y+y\otimes x$ of the subring
$N\subseteq H^*(\EG X^2)$. Such an $\alpha$ satisfies condition
 (a) of Haefliger's Theorem, with $\beta=0$, since
$\Delta_G^*(\alpha)=0$. However it cannot satisfy condition (b),
since $r^*$ is injective on $N$. Hence $[\alpha]:=j^*(\alpha)$ is
non-zero, and it follows that $j^*$ is injective on $N$, so
$H^*(B(X;2))$ contains a subring $
\tilde{N}=\{[\alpha]\mid\alpha\in N\}$ which is isomorphic with
$N$.
\end{proof}

\begin{corollary}
 The subring $\tilde{N}$ of $H^*(B(X;2))$ is contained in
  $\mathrm{Ker}(\pi_G^*)$ where $\pi_G^\ast\co H^\ast(B(X,2)) \to H^\ast(P'X/G)$
  is the map induced by fibration (\ref{map}).
\end{corollary}
\begin{proof}
The generalised diagonal map factorises as \bdm \Delta_G\co
BG\times X\stackrel{\phi}{\to}\EG
PX\stackrel{(1\times_G\pi)}{\longrightarrow}\EG X^2, \edm where
the map $\phi([e],A)=[e,\textrm{const}(A)]$ is a homotopy
equivalence; here $\textrm{const}(A)$ denotes constant path at
$A$. It follows for all $\alpha\in N$ that
$(1\times_G\pi)^*(\alpha)=0$, since $\Delta_G^*(\alpha)=0$.
Applying cohomology to the diagram (\ref{diag1}),
\begin{diagram}[height=2.2em]
H^*(P'X/G) &\lTo &H^*(\EG PX) \\
\uTo<{\pi_G^*} & &\uTo{(1\times_G\pi)^*} \\
H^*(B(X;2)) &\lTo^{j^*} & H^*(\EG X^2)
\end{diagram}
we see that $\pi_G^*[\alpha]=0$ for all $[\alpha]\in\tilde{N}$.

\end{proof}

The following theorem is the main result of this Section:
 \begin{theorem}\label{lower} Let $X$ be a closed smooth manifold.
 Then
\begin{eqnarray}\label{lb1} \TCS(X)\geq\textrm{\rm {cup-length}}(N)+2.
\end{eqnarray}
\end{theorem}
\begin{proof}
This follows by combining the previous Corollary with Theorem
\ref{cohom}.
\end{proof}

\begin{corollary}
For any closed connected smooth manifold $X$ of dimension $\dim
X>0$ one has $\TCS(X)$ $\geq 3$. In particular, $\TCS(S^n)=3$ for
$n\geq 1$.
\end{corollary}
\begin{proof}
The fundamental class $a\in H^n(X)$ gives a non-zero element
$a\otimes 1+1\otimes a$ of $\mathrm{Ker}(\pi_G^*)$, and so
$\TCS(X)\geq 3$. The second statement now follows from Proposition
\ref{sphere}.
\end{proof}

Comparing with the corresponding result for the usual topological
complexity \cite{Far03} we see that $\TC(S^n)=\TCS(S^n)$ for $n$
even while $\TC(S^n)< \TCS(S^n)$ for $n\geq 1$ odd.

\begin{remark}
 The subring $N\subseteq H^*(X)\otimes
H^*(X)$ is contained in the kernel $I\subseteq H^*(X)\otimes
H^*(X)$ of the homomorphism $\Delta^*\co H^*(X\times X)\cong
H^*(X)\otimes H^*(X)\to H^*(X)$ induced by the diagonal map, since
we are working mod 2. This kernel $I$ is precisely the {\em ring
of zero-divisors}, introduced in \cite{Far03}, whose cup-length
plus one provides a lower bound for $\TC(X)$, i.e.
\begin{eqnarray}\label{lb0}\TC(X)\geq \textrm{\rm {cup-length}}(I)+1.\end{eqnarray}
Comparing (\ref{lb1}) with (\ref{lb0}) we see that $\textrm{\rm
{cup-length}}(I)\geq \textrm{\rm {cup-length}}(N)$ however in
(\ref{lb0}) there is an extra 1 which makes estimate (\ref{lb1})
stronger in some cases.
\end{remark}

\section{Cohomology of $P'X$}

In this Section we collect information about cohomology of the
path spaces $P'X$. The results of this Section are not used in
this paper and therefore we will state them without proofs.

Recall that $P'X$ is defined as the subspace of the full path
space $PX$ consisting of all paths $\gamma\co[0,1]\to X$ with
distinct end points $\gamma(0)\not=\gamma(1)$. In other words
$
P'X = PX - LX
$
where $LX$ denotes the space of free loops in $X$, i.e. the space
of all $\gamma\co[0,1]\to X$ with $\gamma(0)=\gamma(1)$. We consider
$P'X$ and $LX$ with the topology induced from $PX$.
\begin{proposition}\label{coh}
Assume that $X$ is a closed connected orientable $n$-dimensional
manifold. Let $\chi$ denote the Euler characteristic $\chi(X)$.
Then
\begin{enumerate}
\item $H^i(P'X;\Z) \simeq H^i(X;\Z)$ for $i<n-1$; \item If
$\chi\not=0$, there is an exact sequence
$$0\to H^{n-1}(X;\Z)\to H^{n-1}(P'X;\Z)\to \tilde H^0(LX;\Z)\to
0;$$ In the case when $\chi=0$ the group on the right should be
replaced by the unreduced zero-dimensional cohomology
$H^0(LX;\Z)$; \item There is an exact sequence
$$0\to \Z/\chi \Z \to H^n(P'X;\Z) \to H^1(LX;\Z)\to 0;$$
 \item One has $H^i(P'X)
\simeq H^{i-n+1}(LX)$ for $i>n$.
\end{enumerate}
\end{proposition}

Below we briefly sketch the proof. One starts with the long exact
sequence
$$\dots \to H^i(PX) \to H^i(P'X) \stackrel \delta\to H^{i+1}(PX,P'X) \to\dots$$
and observes that $H^i(PX)\simeq H^i(X)$ and
$H^{i+1}(PX,P'X)\simeq H^{i+1-n}(LX)$. The second isomorphism is
obtained as a composition of an excision and a Thom isomorphism.
To justify it one notes that $LX$ has an open neighborhood in $PX$
which is homeomorphic to the total space of the rank $n$ vector
bundle which is induced by the map $LX \to X$ (assigning to a loop
its beginning) from the tangent bundle $TX\to X$ over $X$. For
dimensional reasons the connecting homomorphism $$\delta\co
H^i(LX)\to H^{i+n}(X)$$ can be nonzero only for $i=0$; in that
case it coincides with the restriction $H^0(LX)\to H^0(X)$
composed with multiplication by the Euler class of the tangent
bundle.

Note that $H^0(LX;\Z)$ can be identified with the group of all
integer valued functions on the set of conjugacy classes of
$\pi_1(X)$. The action of the natural involution $P'X\to P'X$ on
the cohomology can be described using natural isomorphisms of
Proposition \ref{coh}. For instance we mention that the
isomorphism of statement (4) of Proposition \ref{coh} expresses
the action of the involution on $H^i(P'X)$ for $i>n$ through the
cohomology of the free loop space $LX$ viewed together with the
involution given by reversing directions of loops.

\section{Mid-point maps}\label{midpointsec}
In this Section we study symmetric motion planning algorithms of a
different type which leads to another notion of symmetric
topological complexity of configuration spaces.

Recall that a Symmetric Motion Planner in $X$ is a function $s\co
X\times X\to PX$ which assigns to each pair $(A,B)\in X\times X$ a
path $s(A,B)$ in $X$ from $A$ to $B$, satisfying $s(A,A)(t)=A$ and
$s(B,A)(t)=s(A,B)(1-t)$ for all $t\in I$. Considering just the
mid-points of such paths, if we set
\begin{eqnarray}\sigma
(A,B)=s(A,B)(1/2)\end{eqnarray} we obtain a function $\sigma\co
X\times X\to X$ which satisfies the conditions
$\sigma(A,B)=\sigma(B,A)$ and $\sigma(A,A)=A$.\begin{center}
\resizebox{8cm}{2.5cm}{\includegraphics[55,541][460,670]{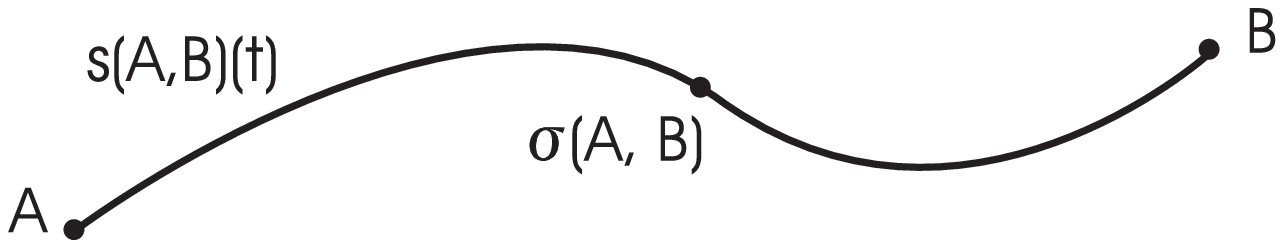}}
\end{center}

A {\em continuous} function $\sigma$ with these properties is
called a ``2-mean'' on $X$. The question of existence of such
means on spaces was considered by B.\ Eckmann, first in 1954
\cite{Eck54} and again 50 years later in relation to the problem
of Social Choice in Economics \cite{Eck04}. For example, any
CW-complex $X$ which admits a 2-mean must be an $H$-space, and
either contractible or of infinite dimension. If we abandon the
requirement $\sigma(A,A)=A$, however, then such maps are not so
rare.
\begin{definition}
A {\em midpoint map} on a topological space $X$ is a continuous map
\begin{eqnarray}\label{midpointmap}
\sigma\co F(X;2)\to X
\end{eqnarray}
satisfying $\sigma(A,B)=\sigma(B,A)$ for all $(A,B)\in F(X;2)$.
\end{definition}
It is clear that such mid-point maps exist for any space $X$ (for
example by setting $\sigma(A,B)=A_0$, where $A_0\in X$ is a base
point) and that they are classified by homotopy classes of continuous
maps from $B(X;2)$ to $X$.
\begin{proposition}
Any midpoint map on $S^n$ with $n$ even is homotopic to a constant
map. For $n$ odd the midpoint maps on $S^n$ are classified by the
integers (the so-called {\em degree} $d\in\Z$).
\end{proposition}
\begin{proof}
The inclusion $S^n\hookrightarrow F(S^n;2)$ is $G$-homotopy
equivalence, where $S^n$ is viewed as a $G$-space with antipodal
action. Hence it induces a homotopy equivalence $\R P^n\to B(S^n;2)$.
The mid-point maps on $S^n$ are thus classified by homotopy
classes of maps from $\R P^n$ to $S^n$, which (by Hopf's Theorem) are
in 1-1 correspondence with $H^n(\R P^n;\mathbb{Z})$. The latter group
is trivial for $n$ even and $\Z$ for $n$ odd, which proves the claim.
\end{proof}
When $n$ is odd there exist mid-point maps on $S^n$ of arbitrary
degree. This is illustrated in the case $n=1$ by considering the map
$\sigma\co F(S^1;2)\to S^1$ given by $\sigma(A,B)=(AB)^d$, where
$d\in\Z$. This uses the structure of an Abelian group on $S^1$.

For any topological space $X$ and a midpoint map
(\ref{midpointmap}) we now define a number $\TCS_\sigma(X)$ which
measures the essential discontinuities of Symmetric Motion
Planners in $X$ whose mid-points are determined by $\sigma$.
Consider the subspace \begin{eqnarray}\label{13.5}
E_\sigma'=\{\gamma\in
 PX\mid\gamma(0)\neq\gamma(1),\gamma(1/2)=\sigma(\gamma(0),\gamma(1))\}\subseteq PX
\end{eqnarray}
consisting of paths with distinct endpoints and mid-point
determined by $\sigma$. We denote by $\pi^\sigma\co E_\sigma'\to
F(X;2)$ the fibration resulting from restricting the endpoint map
(\ref{pi}) to this subspace. This is a $G$-equivariant map of free
$G$-spaces and so induces a fibration
\begin{eqnarray}\label{pisigma}
\pi_G^\sigma:=\pi^\sigma/G\co E_\sigma'/G\to B(X;2).
\end{eqnarray}
\begin{definition}
The {\em Symmetric Topological Complexity of $X$ with
 mid-point map
  $\sigma$} is the number $\TCS_\sigma(X)$ defined as one plus the
Schwarz genus of the fibration $\pi_G^\sigma$.
\end{definition}
\begin{proposition} One has
$
\TCS(X)\leq \TCS_\sigma(X).
$
\end{proposition}
\begin{proof}
The fibrations $\pi_G$ and $\pi_G^\sigma$ are related as shown in
the diagram
\begin{eqnarray*}
\begin{array}{ccc}
 E_\sigma'/G  &\stackrel{\varphi}\longrightarrow  & P'X/G \\
 \qquad\qquad\qquad \searrow{\pi_G^\sigma}\quad & & \swarrow{\pi_G}\qquad\qquad\qquad \\
 &  B(X;2)  &\end{array}
\end{eqnarray*}
where the map $\varphi$ is induced by the inclusion
$E_\sigma'\subseteq P'X$. Any section
$s\co U\to E_\sigma'/G$ of $\pi_G^\sigma$ on an open set
$U\subseteq B(X;2)$ therefore gives rise to a section
$\varphi\circ s$ of $\pi_G$ on $U$, and the conclusion follows.
\end{proof}

The proof of Proposition \ref{dim} works with $\TCS_\sigma(X)$
replacing $\TCS(X)$ and gives:
\begin{proposition}\label{dim1} For
any $X$ we have $\TCS_\sigma (X)\leq 2\mathrm{dim}(X)+2$. If $X$
is a closed smooth manifold then $\TCS(X)\leq 2\mathrm{dim}(X)+1$.
\end{proposition}

The main result of this paper concerning $\TCS_\sigma(X)$ can be
stated as follows:

\begin{theorem}\label{cuplength}
Let $X$ be a smooth closed aspherical manifold. Then
\begin{eqnarray}
\TCS_\sigma(X) \geq 2 \cl(X) +1
\end{eqnarray}
where $\sigma\co F(X;2)\to X$ is the constant mid-point map and
$\cl(X)$ denotes the largest integer $k$ such that there exist $k$
cohomology classes $u_1, \dots, u_k\in H^\ast(X;\Z_2)$ of positive
degree whose cup-product is non-zero, $u_1\dots u_k\not=0$.
\end{theorem}

Recall that a path-connected topological space $X$ is said to be
aspherical if $\pi_i(X)=0$ for all $i>1$.

Theorem \ref{cuplength} will be proven in the final Section
\ref{secaspherical}. We first need to sharpen our tools by
introducing a new variant on category weight. This will be done in
Section \ref{secsec}.

We conclude this Section by two examples:

\begin{example}{\rm
Consider the closed orientable surface $\Sigma_g$ of genus $g\geq
1$. This has cup-length $\cl(\Sigma_g)=2$, and so
$\TCS_\sigma(\Sigma_g)\geq 5$ by Theorem \ref{cuplength}.
Proposition \ref{dim1} gives $\TCS_\sigma(X)\leq 2n+1$ whenever
$X$ is a closed $n$-manifold. Hence $\TCS_\sigma(\Sigma_g)=5$.
This agrees with the usual Topological Complexity except in the
case of the torus $T^2=\Sigma_1$ which has $\TC(T^2)=3$, see
\cite{Far03}.}
\end{example}
\begin{example}{\rm
Consider the $n$-dimensional torus $T^n$, the Cartesian product of
$n$ copies of $S^1$.  This models the configuration space of a
planar robot arm with $n$ revolving joints. It has $\cl(T^n)=n$,
and Theorem \ref{cuplength} implies that $\TCS_\sigma(T^n)=2n+1$.
The usual topological complexity of the torus is $\TC(T^n)=n+1$. }
\end{example}

The last example shows that $\TCS_\sigma(X)$ can be much larger
than $\TC(X)$.

\section{Sectional category weight}\label{secsec}
The usual cohomological lower bound for the Lusternik-Schnirelmann
category of a space $X$ states that if $u_1,\ldots ,u_k\in
H^\ast(X)$ are non-zero cohomology classes of positive degree such
that their cup product $u=u_1\cdots u_k\neq 0$ is nonzero then
$\mathrm{cat}(X)> k$. E.\ Fadell and S.\ Husseini \cite{FH92}
improved this estimate by assigning to each cohomology class $u$
an integer weight, denoted $\mathrm{cwgt}(u)$ and called its {\em
category weight}, such that if $0\neq u=u_1\cdots u_k$ then \bdm
\mathrm{cat}(X)>\mathrm{cwgt}(u)=\mathrm{cwgt}(u_1\cdots
u_k)\geq\sum_i \mathrm{cwgt}(u_i). \edm Note that to improve on
the cup-length estimate one must find indecomposable classes of
category weight at least 2, which Fadell and Husseini did using
Steenrod operations \cite{FH92}. This notion was developed further
by Y.\ B.\ Rudyak, who gave a homotopy invariant version (known as
{\em strict category weight}) and proved that non-trivial Massey
products also have weight at least 2 \cite{Rud99}.

Our aim here is to improve on Schwarz's original cohomological lower
bound for the genus of a fibration $p\co E\to B$
(given in Theorem \ref{cohom}) by extending the ideas of Fadell and Husseini
mentioned above. To do this we note that each cohomology class
$\xi\in H^*(B)$ can be assigned an integer weight {\em with respect to
  $p$}, called its {\em sectional category weight}, and defined as follows.

If $f\co X\to B$ is a continuous map, we denote by $f^*p\co
f^*(E)\to X$ the pull-back fibration; the symbol $g(f^\ast p)$
denotes the Schwarz genus of $f^\ast p$.

\begin{definition}
 Let $\xi\in H^*(B)$ be a cohomology
 class. We define the sectional category weight of $\xi$ with respect
 to $p$, denoted $\wgt_p(\xi)$, to be the largest integer $k$ such that
 for any continuous map $f\co X\to B$ with $g(f^*p)\leq k$ we have
 $f^*(\xi)=0$.
\end{definition}
Note that $\wgt_p(\xi)\geq 0$, since any fibration over a non-empty base
has genus at least
$1$, and so the condition $f^*(\xi)=0$ whenever $g(f^*p)\leq 0$ is
vacuously satisfied. We set $\wgt_p(\xi)=\infty$ when $\xi=0$.
\begin{proposition}\label{weight1}
A class $\xi\in H^*(B)$ has $\wgt_p(\xi)\geq 1$ if and only if
$p^*(\xi)=0\in H^*(E)$.
\end{proposition}
\begin{proof}
First suppose $p^*(\xi)=0$, and let $f\co X\to B$ be a map such that
$g(f^*p)\leq 1$, so that $f^*p$ has a section. From the pullback
diagram
\begin{diagram}[height=2em]
f^*(E) &\rTo^{\overline{f}} & E \\
\dTo<{f^*p} & &\dTo>{p} \\
X & \rTo^{f} & B
\end{diagram}
we see that $(f^*p)^*f^*(\xi)=\overline{f}^*p^*(\xi)=0$, which implies
$f^*(\xi)=0$ since $(f^*p)^*$ is injective. Hence $\wgt_p(\xi)\geq 1$.

Conversely suppose that $\wgt_p(\xi)\geq 1$, and consider the pull-back
diagram
\begin{diagram}[height=2em]
p^*(E) &\rTo & E \\
\dTo<{p^*p} & & \dTo>{p} \\
E &\rTo^{p} & B.
\end{diagram}
Clearly the fibration $p^*p$ has a section given by the diagonal map,
and so $p^*(\xi)=0$ by the
definition of sectional category weight.
\end{proof}

\begin{proposition}
For any non-zero class $\xi\in H^*(B)$,
\begin{eqnarray}
\wgt_p(\xi)< g(p).
\end{eqnarray}
\end{proposition}
\begin{proof}
Suppose the converse is true, and we have $\xi$ such that
$\wgt_p(\xi)\geq g(p)=k$, say. The identity map $\mathbf{1}\co B\to B$
has $g(\mathbf{1}^*p)=g(p)=k$, and so $\mathbf{1}^*\xi =\xi$ must be zero.
\end{proof}
\begin{proposition}
For the cup product of classes $\xi_1,\ldots ,\xi_l\in H^*(B)$ we have
\begin{eqnarray}
\wgt_p(\xi_1\cdots\xi_l)\geq\sum_{i=1}^{l}\wgt_p(\xi_i).
\end{eqnarray}
\end{proposition}
\begin{proof}
Of main interest is the case when $\xi=\xi_1\cdots\xi_l\neq 0$ (the other
case is true by convention). Letting $k_i=\wgt_p(\xi_i)$, we must show
that $\wgt_p(\xi)\geq k=\sum k_i$. So suppose that $f\co X\to B$ is a
continuous map with $g(f^*p)\leq k$. We may find a covering
$\Omega=\{ U_1,\ldots ,U_k\}$ of $X$ by open sets, above each of which the
fibration $f^*p$ has a section. Partition the cover $\Omega$ into
$l$ families $\Omega_1,\ldots ,\Omega_l$ such that $\Omega_i$ consists
of $k_i$ open sets. We now define
\bdm
A_i=\bigcup_{U_j\in\Omega_i}U_j.
\edm
Note that $X=\bigcup_i A_i$, and $f^*(\xi_i)|_{A_i}=0$ for
$i=1,\ldots ,n$ since $\wgt_p(\xi_i)=k_i$. A standard argument now gives that
$
f^*(\xi)=f^*(\xi_1)\cdots f^*(\xi_l)=0,
$
as required.
\end{proof}
The last three Propositions together imply the following
sharpened version of Theorem \ref{cohom}.
\begin{theorem} \label{Sweight}
Let $p\co E\to B$ be a fibration. If $\xi_1,\ldots ,\xi_l\in H^*(B)$ are positive dimensional cohomology classes whose
product is non-zero, then
\begin{eqnarray}
g(p)>\sum_{i=1}^l \wgt_p(\xi_i).
\end{eqnarray}
\end{theorem}
Hence we may improve on the lower bound given by Theorem \ref{cohom} by
finding indecomposable elements of sectional category weight at least 2. Finding such elements is greatly facilitated by
the next Proposition regarding fibrewise joins.

Recall that the {\em fibrewise join} of two fibrations $p_1\co E_1\to
B$ and $p_2\co E_2\to B$ over the same base is a certain fibration
$p_1\ast p_2\co E_1\ast_B E_2\to B$, whose fibre has the homotopy type
of the join $F_1\ast F_2$ (see \cite{Sch66} or \cite{Ark99}, for
example). The total space $E_1\ast_B E_2$ may be described as the subspace
\bdm
\{(e_1,e_2,t)\in E_1\times E_2\times I\mid
p_1(e_1)=p_2(e_2)\}\subseteq E_1\times E_2\times I
\edm
modulo the relations $(e_1,e_2,0)\sim (e_1',e_2,0)$ and
$(e_1,e_2,1)\sim (e_1,e_2',1)$ for all $e_1,e_1'\in E_1$ and
$e_2,e_2'\in E_2$. The projection $p_1*p_2$ of this fibration is given by
$p_1\ast p_2([e_1,e_2,t])=p_1(e_1)=p_2(e_2)$.

One may also define the fibrewise join of an arbitrary number of
fibrations. For a fibration $p\co E\to B$, denote by $p(k)\co
E(k)\to B$ the $k$-fold fibrewise join of $p\co E\to B$ with
itself. A. Schwarz (\cite{Sch66}, Theorem 3) proved that {\em
$g(p)\leq k$ if and only if $p(k)$ has a continuous section}.
\begin{proposition}\label{join}
If $p(k)^*(\xi)=0$ then $\wgt_p(\xi)\geq k$.
\end{proposition}
\begin{proof}
Let $f\co X\to B$ be such that
$f^*p\co f^*(E)\to X$ has genus not more than $k$. The fibrations $f^*p(k)$ and $(f^*p)(k)$ are homeomorphic. Hence we have the following diagram,
\begin{diagram}[height=2em]
f^*(E)(k) &\rTo & E(k) \\
\dTo<{(f^*p)(k)} & & \dTo>{p(k)} \\
 X & \rTo^{f} & B
\end{diagram}
in which the map $(f^*p)(k)$ admits a section and hence induces a
monomorphism in cohomology. If $p(k)^*(\xi)=0$ then
$(f^*p)(k)^*f^*(\xi)=0$ and hence $f^*(\xi)=0$.
\end{proof}
To conclude this Section we describe the homotopy
invariance of sectional category weight.
\begin{lemma}\label{htpy}
Suppose $p_1\co E_1\to B$ and $p_2\co E_2\to B$ are fibre homotopy
equivalent fibrations, and $\xi\in H^*(B)$. Then
$
\wgt_{p_1}(\xi)=\wgt_{p_2}(\xi).
$
\end{lemma}
\begin{proof}
Given any continuous map $f\co X\to B$, the pullback fibrations
$f^*p_1$ and $f^*p_2$ are fibre homotopy equivalent. Since fibre
homotopy equivalent fibrations have the same genus (this follows
immediately from Proposition 6 of Schwarz \cite{Sch66}), this means that
$g(f^*p_1)=g(f^*p_2)$ for all $f$, and the conclusion follows.
\end{proof}
\begin{proposition}\label{pullback}
Let $p\co E\to B$ be a fibration, $\xi \in H^\ast(B)$ a cohomology
class and $g\co A\to B$ a continuous map. Then, $\wgt_{g^*p}g^*(\xi)\geq\wgt_p(\xi). $ If
$g$ is a homotopy equivalence, then
$\wgt_{g^*p}g^*(\xi)=\wgt_p(\xi)$.
\end{proposition}
\begin{proof}
The first statement is immediate from the definition. Suppose
$k\co B\to A$ is a homotopy inverse for $g$. Then $g\circ k\simeq
1$ implies $k^*g^*(\xi)=\xi$ and $k^*g^*p$ is fibre homotopy
equivalent to $p$, and so by Lemma \ref{htpy} \bdm
\wgt_p(\xi)=\wgt_{k^*g^*p}k^*g^*(\xi)\geq\wgt_{g^*p}g^*(\xi)\geq\wgt_p(\xi).
\edm
\end{proof}

\section{Proof of Theorem \ref{cuplength}}\label{secaspherical}

Recall that the number $\TCS_\sigma(X)$ is defined to be one plus
the genus of
 the fibration (\ref{pisigma}).  The space
\bdm E_\sigma =\{\gamma\co I\to X\mid \gamma(1/2)=A_0\}\subseteq
PX \edm of paths in $X$ with mid-point $A_0$ admits an involution
$\gamma\mapsto\overline{\gamma}$ where
$\overline{\gamma}(t)=\gamma(1-t)$, and it contains $E_\sigma'$
(given by (\ref{13.5})) as a free $G$-invariant subspace (the
subspace consisting of such paths with distinct endpoints). The
endpoint map $q\co E_\sigma\to X^2$ is a $G$-invariant fibration,
of which $\pi^\sigma\co E_\sigma'\to F(X;2)$ is a restriction. We
obtain a diagram of homotopy orbit spaces

\begin{eqnarray}\label{EGBG}
\begin{array}{ccccc}
E_\sigma'/G & \simeq & EG\times_G E_\sigma' & \to & \EG E_\sigma
\\ \\
\quad \downarrow {\pi_G^\sigma} & &\qquad\qquad \downarrow
{1\times_G\pi^\sigma} &
&\qquad\qquad\downarrow {1\times_G q=p} \\ \\
B(X;2) &\simeq & \EG F(X;2) &\stackrel j\to & \EG X^2.
\end{array}
\end{eqnarray}
By Theorem \ref{Steen} and the subsequent paragraphs, every
cohomology class $u\in H^k(X)$ determines cohomology classes
\begin{eqnarray}
\alpha_u = 1\otimes u\otimes u\in H^{2k}(EG\times_G X^2 ), \quad
\beta_u=j^\ast(\alpha_u)\in H^{2k}(B(X,2)).
\end{eqnarray}
The main ingredient of the proof of Theorem \ref{cuplength}
consists of the following statement:

\begin{theorem}\label{atleast} For any $u\in H^k(X)$ with
$k>0$, the cohomology class $\alpha_u\in
H^{2k}(EG\times_G X^2)$ has sectional category
weight at least $2$ with respect to the fibration $p\co
EG\times_G E_\sigma \to EG\times_G X^2$. Hence, the cohomology
class $\beta_u\in H^{2k}(B(X,2))$ has sectional category weight at
least $2$ with respect to the fibration $\pi^\sigma_G\co E'_G/G \to
B(X,2)$.
\end{theorem}
The proof of Theorem \ref{atleast} rests on the following two
lemmas, the first of which concerns
 the fibrewise join
 $$q\ast q\co E_\sigma
  *_{X^2}E_\sigma\to X^2 $$
  of two copies of the fibration
  $q\co E_\sigma\to X^2 $. Let $\Omega X$ denote the space of
  loops in $X$ based at the base-point $A_0\in X$. We think of points
  of $\Omega X$ as maps $\omega\co I\to X$ such that
  $\omega(0)=\omega(1)=A_0$. Let $S$ denote the unreduced
  suspension functor. There is a map
\bdm h\co S(\Omega X\times \Omega X)\to X^2  \edm defined by
setting $h[\omega_1,\omega_2,t]=(\omega_1(1/2),\omega_2(1/2))$.
\begin{lemma}\label{loops}
There is a homotopy equivalence $ \tau\co E_\sigma*_{X^2}E_\sigma
\to S(\Omega X\times\Omega X) $ which satisfies $h\circ\tau =q*q$.
Furthermore, $\tau$ is $G$-equivariant with respect to the
diagonal action of $G$ on $E_\sigma*_{X^2}E_\sigma$ and the action
on $S(\Omega X\times\Omega X)$ which swaps the loop factors.
\end{lemma}
\begin{proof}[Proof of Lemma \ref{loops}](Compare Theorem 21
of \cite{Sch66}.) Points of $E_\sigma*_{X^2}E_\sigma$ are
equivalence classes $[\gamma_1,\gamma_2,t]$ where $\gamma_1$,
$\gamma_2$ are paths in $X$ with the same initial and final points
and mid-point $A_0$ (see figure below). We set
$\tau[\gamma_1,\gamma_2,t]=[\omega_1,\omega_2,t]$, where

$$ \omega_1(s)=\left\{ \begin{array}{ll}
         \gamma_2(1/2-s) &\textrm{if $s\in[0,1/2)$},\\
         \gamma_1(s-1/2) &\textrm{if $s\in[1/2,1]$},
\end{array}\right., \
 \omega_2(s)=\left\{ \begin{array}{ll}
         \gamma_2(1/2+s) &\textrm{if $s\in[0,1/2)$},\\
         \gamma_1(3/2-s) &\textrm{if $s\in[1/2,1]$}.
\end{array}\right.
$$
We then have $\tau
[\overline{\gamma}_1,\overline{\gamma}_2,t]=[\omega_2,\omega_1,t]$
and $h\circ\tau=q*q$ by design. \begin{center}
\resizebox{9cm}{4cm}{\includegraphics[34,447][567,712]{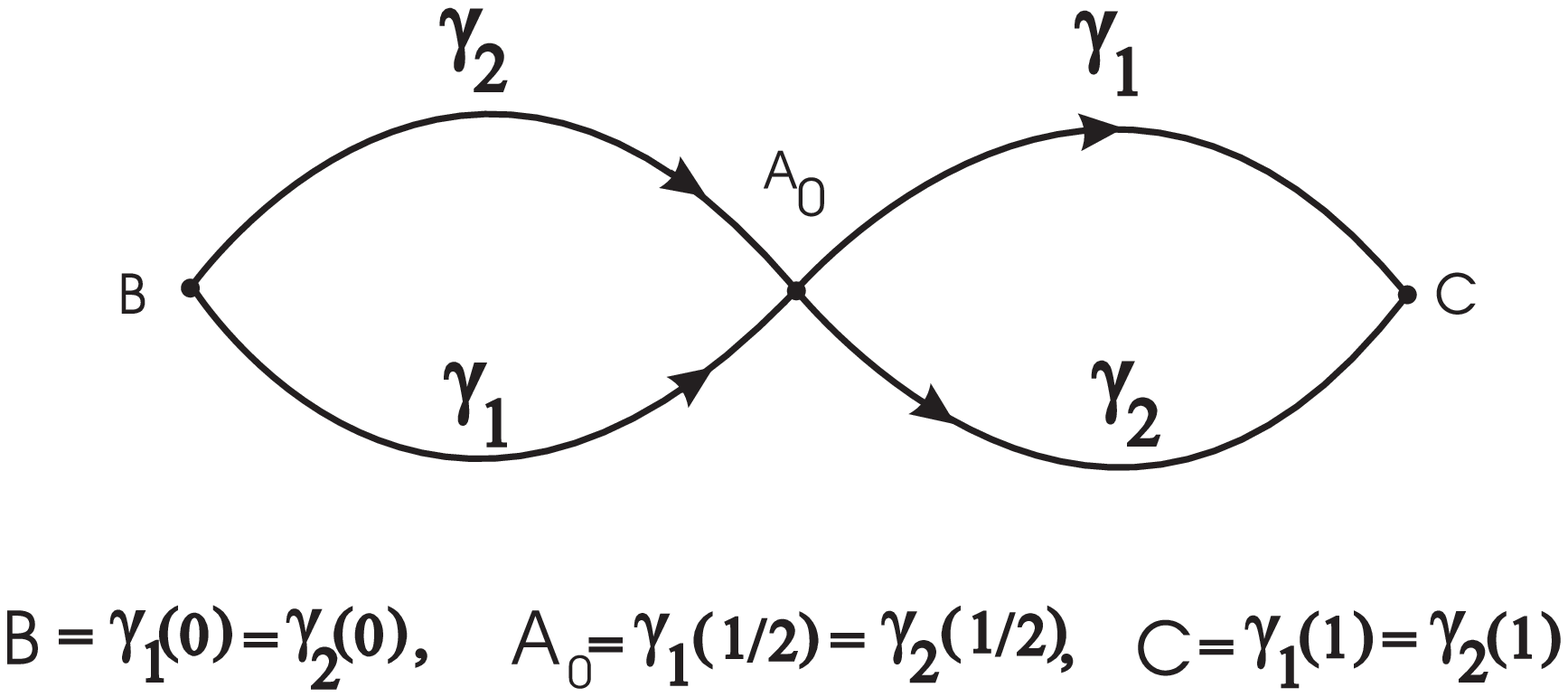}}
\end{center}
It remains to check that $\tau$ is a homotopy equivalence. Denote
by $W_0$ (respectively $W_1$) the closed subspace of
$E_\sigma*_{X^2 }E_\sigma$ consisting of points of the form
$[\gamma_1,\gamma_2,0]$ (resp.\ $[\gamma_1,\gamma_2,1]$). Both of
these subspaces are homeomorphic to $E_\sigma$, hence are
contractible within themselves. It is not difficult to construct a
deformation of the space $E_\sigma*_{X\times
  X}E_\sigma$ which contracts each to a point within itself. The proof
is completed by noting that, up to homeomorphism, $\tau$ is the
quotient map collapsing $W_0$ and $W_1$ to the vertices of the
suspension, and so is a homotopy equivalence by Lemma 7.1.5 of
Spanier \cite{Spa66}.
\end{proof}

Our next result generalises to the
equivariant setting Proposition 1 of Schwarz \cite{Sch66}
concerning the join of two fibrations associated with a principle
fibration.

 Let $E$ be a free $G$-space and let $q_i\co Y_i\to X$,
$i=1,2$ be two equivariant maps of $G$-spaces which are fibrations
over the common base space $X$. The group $G$ acts diagonally on
the total space $Y_1*_XY_2$ of the fibrewise join $q_1*q_2$.

\begin{lemma}\label{homeo} The following two fibrations over $E\times_G X$ are
homeomorphic: $$ 1\times_G(q_1*q_2)\co E\times_G(Y_1*_XY_2)\to
E\times_G X, $$ and $$ (1\times_Gq_1)*(1\times_Gq_2)\co (E\times_G
Y_1)*_{E\times_GX}(E\times_G Y_2)\to E\times_G X. $$
\end{lemma}
\begin{proof}[Proof of Lemma \ref{homeo}]
We define a map \begin{eqnarray}\label{mu} \mu\co
E\times_G(Y_1*_XY_2)\to (E\times_G Y_1)*_{E\times_GX}(E\times_G
Y_2)\end{eqnarray} and show that it is a homeomorphism. A point in
the domain of (\ref{mu}) is an equivalence class
$[e,[y_1,y_2,t]]$, where $e\in E$, $y_i\in Y_i$ and $t\in I$ such
that $q_1(y_1)=q_2(y_2)\in X$. For any $g\in G$ we have
$[e,[y_1,y_2,t]]=[ge,[gy_1,gy_2,t]]$. A point in the range is a
class $[[e_1,y_1],[e_2,y_2],t]$ where $e_1,e_2\in E$, $y_i\in Y_i$
and $t\in I$ such that $[e_1,q_1(y_1)]=[e_2,q_2(y_2)]\in
E\times_GX$. We set \bdm \mu[e,[y_1,y_2,t]]=[[e,y_1],[e,y_2],t].
\edm This map clearly is well-defined and is continuous. We will
show that $\mu$ is onto; the proof that it's 1-1 runs similarly.

A point $x=[[e_1,y_1],[e_2,y_2],t]$ in the range has
$[e_1,q_1(y_1)]=[e_2,q_2(y_2)]\in E\times_GX$. Hence there is a
unique $g\in G$ with $e_1=ge_2$ and $q_1(y_1)=gq_2(y_2)$.
Therefore,
\begin{eqnarray*}
x & = & [[e_1,y_1],[e_1,gy_1],t]\\
 & = & \mu[e_1,[y_1,gy_2,t]].
\end{eqnarray*}
(Note that in this way one may describe a continuous inverse for
$\mu$).

Hence $\mu$ is a homeomorphism. Verifying that $\mu$ is fibre
preserving is trivial.
\end{proof}

\begin{proof}[Proof of Theorem \ref{atleast}] We will show that
$ (p\ast p)^\ast(\alpha_u)=0. $ The first statement of Theorem
\ref{atleast} then follows from Proposition \ref{join}, and the
second statement follows from the first using Proposition
\ref{pullback} since $\pi^\sigma_G$ is the pullback fibration
$j^\ast p$ where $j\co B(X,2)\to EG\times_G X^2$ is the inclusion.
Applying Lemma \ref{homeo} we find that $(p\ast
p)^\ast(\alpha_u)=0$ if and only if $(1\times_G(q\ast
q))^\ast(\alpha_u)=0$.

Lemma \ref{loops} gives a homotopy equivalence \bdm \tau\co
E_\sigma*_{X^2}E_\sigma \simeq S(\Omega X\times\Omega X). \edm
Since $X$ is an aspherical manifold, it is a $K(\pi,1)$ with
discrete fundamental group. Thus $\Omega
X\times\Omega X$ has the homotopy type of a discrete set of
points, and $S(\Omega X\times\Omega X)$ of a wedge of circles.
Therefore the map $(q\ast q)^*\co H^{2k}(X^2 )\to
H^{2k}(E_\sigma*_{X^2 }E_\sigma)$ takes values in a zero
group when $k>0$; in particular, $(q\ast q)^*(u\otimes u)=0$.

The map $q\ast q\co E_\sigma*_{X^2 }E_\sigma\to X^2 $ is
$G$-equivariant, and so gives a morphism of the associated
Cartan-Leray spectral sequences (see \cite{Bor60}, Chapter IV,
Section 3),
\begin{diagram}[height=2.1em]
H^p(G,H^q(E_\sigma*_{X^2 }E_\sigma)) &\implies &
H^{p+q}(\EG (E_\sigma*_{X^2 }E_\sigma))\\
\uTo<{ (q\ast q)^*} & &\uTo>{(1\times_G (q\ast q))^*}\\
H^p(G,H^q(X^2 ))&\implies & H^{p+q}(\EG X^2 ).
\end{diagram}
The lower spectral sequence collapses at the $E_2$ term, by
Theorem \ref{Steen}. Hence we may view $\alpha_u=1\otimes u\otimes
u\in H^{2k}(EG\times_G X^2)$ as an element of $H^0(G,H^{2k}(X^2 ))$, which maps to zero
under the map on coefficients induced by $q\ast q$. The claim follows
by naturality.
\end{proof}

\begin{proof}[Proof of Theorem \ref{cuplength}]
Suppose ${\rm {cl}}(X)=m$, and let $u_1,\ldots,u_m\in H^*(X)$ be
positive dimensional classes with non-zero product. Let
$\alpha_i$ and $\beta_i$ denote the elements $\alpha_{u_i}$ and
$\beta_{u_i}$. Note that whilst the product $\alpha_1\cdots
\alpha_m\in H^*(EG \times_G X^2)$ is non-zero, the product
$$\beta_1\cdots\beta_m =j^*(\alpha_1\cdots\alpha_m)\in H^*(B(X;2))$$
may vanish. However, letting $\gamma_m=[u_m\otimes 1+1\otimes
u_m]\in H^*(B(X;2))$ we find that the product
$$\beta_1\cdots\beta_{m-1}\gamma_m=[u_1\cdots u_m\otimes u_1\cdots
u_{m-1}+u_1\cdots u_{m-1}\otimes u_1\cdots u_m]$$
is non-zero, by Corollary \ref{Nsurvives}.

We now observe that $\wgt_{\pi^\sigma_G}\gamma_m\geq 1$. Indeed, since $EG\times
E_\sigma$ is a contractible space with free $G$-action, there is a
homotopy equivalence $EG\times_G E_\sigma\simeq BG$. Now applying
cohomology to the diagram (\ref{EGBG}) we obtain
\begin{diagram}[height=2em]
H^*(E_\sigma'/G) &\lTo & H^*(BG) \\
\uTo<{(\pi_G^\sigma)^*} & &\uTo>{\Theta^*} \\
H^*(B(X;2)) & \lTo^{j^*} & H^*(\EG X^2),
\end{diagram}
where $\Theta\co BG\to\EG X^2$ is the map $[e]\mapsto
[e,A_0,A_0]$. Now $\Theta^*(u_m\otimes 1+1\otimes u_m)=0$, which
implies $(\pi_G^\sigma)^*(\gamma_m)=0$. This proves the claim in
view of Proposition \ref{weight1}.

The proof of Theorem \ref{cuplength} is completed by applying Theorem
\ref{Sweight} to the product $\beta_1\cdots\beta_{m-1}\gamma_m$, since
$$g(\pi_G^\sigma)>\sum_{i=1}^{m-1}\wgt_{\pi_G^\sigma}\beta_i
+\wgt_{\pi_G^\sigma}\gamma_m\geq 2(m-1)+1=2m-1.$$
 Hence $g(\pi_G^\sigma)\geq 2m$ and so $\TCS_\sigma(X)\geq 2m+1$, as
stated.
\end{proof}

\end{document}